\documentclass[10pt, a4paper]{article}
\author{Frank Schuhmacher\footnote{Supported by:
Doktorandenstipendium des Deutschen Akademischen Austauschdienstes im Rahmen
des gemeinsamen Hochschulsonderprogramms III des Bundes und der
L\"ander}}

\title{Deformation of Singularities via $L_\infty$-Algebras}

\usepackage{mathenv}
\usepackage[intlimits]{amsmath}
\usepackage{amssymb}
\usepackage{amsthm} 
\usepackage{fancyheadings}
\usepackage{amscd}
\usepackage[all]{xypic}
\usepackage{makeidx}

\newcounter{punkt}

\theoremstyle{definition}
\newtheorem{defi}{Definition}[section]

\newtheorem{prop}[defi]{Proposition}

\theoremstyle{theorem}
\newtheorem{satz}[defi]{Theorem}

\newtheorem{lemma}[defi]{Lemma}
\newtheorem{kor}[defi]{Corollary}

\newcommand{\nach}{\longrightarrow}

\newcommand{\isom}{\cong}

\newcommand{\ZZ}{\mathbb{Z}}

\newcommand{\CC}{\mathbb{C}}

\newcommand{\Oh}{\mathcal{O}}

\newcommand{\ppp}{\cdot\ldots\cdot}

\newcommand{\Mod}{\mathfrak{Mod}}

\newcommand{\sX}{\mathcal{X}}
\newcommand{\im}{\operatorname{im}}
\renewcommand{\Im}{\operatorname{Im}}

\newcommand{\sC}{\mathcal{C}}

\newcommand{\Y}{\mathcal{Y}}

\newcommand{\m}{\mathfrak{m}}

\newcommand{\Hom}{\operatorname{Hom}}

\renewcommand{\dim}{\operatorname{dim}}

\newcommand{\gr}{\operatorname{gr}}

\newcommand{\lz}{\hfill\newline}

\newcommand{\ot}{\otimes}

\newcommand{\Coder}{\operatorname{Coder}}

\newcommand{\dgmanf}{\text{\texttt{DG-Manf}}}

\newcommand{\Der}{\operatorname{Der}}

\newcommand{\Diff}{\operatorname{Diff}}
\newcommand{\Codiff}{\operatorname{Codiff}}

\newcommand{\bR}{\bar{R}}
\newcommand{\an}{\mathfrak{An}}
\newcommand{\anf}{\mathfrak{Anf}}

\newcommand{\obj}{\operatorname{Ob}}
\newcommand{\sD}{\mathcal{D}}
\newcommand{\sY}{\mathcal{Y}}
\newcommand{\tsX}{\tilde{\mathcal{X}}}

\begin{document}

\maketitle
\begin{abstract}
This is an addendum to the paper
``Deformation of $L_\infty$-Algebras'' \cite{Schuh1}.
We explain in which way the deformation theory of $L_\infty$-algebras
extends the deformation theory of singularities.
We show that the construction of semi-universal deformations
of $L_\infty$-algebras gives explicit formal semi-universal
deformations of isolated singularities.  
\end{abstract}

\section*{Introduction}
In this paper, we apply the following general idea
for the construction of moduli spaces
to isolated singularities:
Take the differential graded Lie algebra $L$ describing a deformation
problem (for isolated singularities, this is the tangent complex)
and find a minimal representative $M$ of $L$ in the class of formal
$L_\infty$-algebras (see \cite{Schuh1}). In geometric
terms, $M$ is a formal DG-manifold, containing the moduli space 
as analytic substructure. This general concept is also
sketched in \cite{Merk1}.\\ 

We define a functor $F$ from the category of complex
analytic space germs to the localization of the
category of $L_\infty$-algebras by $L_\infty$-equivalence.
For a singularity $X$,
we take the semi-universal $L_\infty$-deformation $(V,Q^V)$ of $F(X)$ 
constructed in \cite{Schuh1}. For isolated singularities,
the components $V^i$ are of finite dimension.
The restriction of
the vectorfield $Q^V$ defines a formal map (Kuranishi-map)
$V^0\nach V^1$ whose zero locus gives the formal
moduli space.  

\section{Definitions and reminders}

In the whole paper, we work over a ground field
$k$ of characteristic zero.

Denote the category of formal (resp. convergent) complex
analytic space germs by $\anf$ (resp. $\an$).
Denote the category of isomorphism classes
of formal DG manifolds by $\dgmanf$.
We use the following superscripts to denote full
subcategories of $\dgmanf$:\\
L (``local''): the subcategory of all $(M,Q^M)$ in $\dgmanf$
such that $Q^M_0=0$;\\
M (``minimal''): the subcategory of all $(M,Q^M)$ in $\dgmanf^L$
such that $Q^M_1=0$;\\
G (``g-finite''): the subcategory of all $(M,Q^M)$ in $\dgmanf^L$ such that 
$H(M,Q^M_1)$ is g-finite.\\

We call a morphism $f=(f_n)_{n\geq 1}$
in $\dgmanf^L$ \textbf{weak equivalence},
if the morphism $f_1$ of DG vectorspaces is a quasi-isomorphism,
i.e. if the corresponding morphism of $L_\infty$-algebras is
an $L_\infty$-equivalence. Recall that by Theorem 4.4 and
Lemma 4.5 of \cite{Kont},
weak equivalences define an equivalence relation
in $\dgmanf^L$ and that in each equivalence class, there is
a uniquely defined \textbf{minimal model}, i.e. an object belonging to
$\dgmanf^M$. 

\begin{prop}
We can localize the category $\dgmanf^L$ by weak equivalences ($\approx$).
The quotient $\dgmanf^L/\approx$ is equivalent to
the category $\dgmanf^M$ and the localization functor
assigns to each object of $\dgmanf^L$ its minimal
model.
\end{prop}
\begin{proof}
This follows directly by
Corollary 2.5.7 of \cite{Merk1}.
\end{proof}

\section{The functors $F$ and $V$}\label{functors}

In this section we explain how to represent (formal) singularities
by formal DG manifolds.\\

Let $\sC$ be the category of formal analytic algebras,
$A\in\obj(\sC)$ and $R=(R,s)$ a \textbf{resolvent} of $A$ over $k$,
i.e. a g-finite free DG-algebra in $\gr(\sC)$ such that
$H^0(R,s)\isom A$ and $H^j(R,s)=0$, for $j<0$. 
For $l\geq 0$, let $I_l$ be an index set containing one index
for each free algebra generator
of $R$ of degree $-l$. Consider the disjoint union $I$ 
of all $I_l$ as graded set such that $g(i)=l$, for $i\in I_l$.
Fix an ordering on $I$, subject to
the condition $i<j$, if $g(i)<g(j)$.

Thus, as graded algebra, $R=k[[X^0]][X^-]$, 
where $X^0=\{x_i|\;i\in I,\,g(i)=0\}$ 
and $X^-=\{x_i|\;i\in I, g(i)\geq 1\}$ are sets of free algebra
generators with $g(x_i)=-g(i)$.

Set $M:=\coprod_{i\in I}ke_i$ to be the free, graded $k$-vectorspace with
base $\{e_i:i\in I\}$, where $g(e_i)=g(i)$. 
Consider $S(M)=\coprod_{n\geq 0}M^{\odot n}$ in the usual way
as graded coalgebra (see Section 1.1 of \cite{Schuh1}).
Set 
$$S(M)^\ast:=\Hom_{k-\Mod}(S(M),k)=\prod_{j\geq 0}
\Hom_{k-Mod}(M^{\odot j},k).$$

We identify  products $x_{i_1}\cdot\ldots \cdot x_{i_l}$ in $R$
with the maps $M^{\odot l}\nach k$, defined by
$e_{i_1}\ppp e_{i_l}\mapsto 1$ and
$e_{j_1}\ppp e_{j_l}\mapsto 0$ for $\{j_1,\ldots ,j_l\}\neq\{i_1,\ldots ,i_l\}$.
Especially, we identify each constant $\lambda\in k$ with
the map $k\nach k$, sending $1$ to $\lambda$.
We have $$R^j=\prod_{n\geq 0}\Hom^j(M^{\odot n},k)$$ and
$R=\coprod_{j\leq 0}R^j$. The differential $s$ of $R$ extends
naturally to $\bR:=\prod_{j\leq 0}R^j$.
As complexes, $R$ and $\bR$ are identical, but not as graded modules.
We identify $\bR=S(M)^\ast$.
Set 
\begin{equation*}
\Der(R):=\coprod_{i\in\ZZ}\Der^i(R,R)
\quad\text{ and }\quad
\Coder(S(M)):=\coprod_{i\in\ZZ}\Coder^i(S(M),S(M)).
\end{equation*}
Denote $\Diff(R)$ (resp. $\Codiff(S(M))$) the submodule
of differentials (resp. codifferentials).
The following proposition explains why, for a formal DG manifold
$W$, the complex $\Coder(S(W),S(W))$ is called tangent complex
of $W$.

\begin{prop}\label{tan}
Take $R$ and $M$ as above.
The natural map 
\begin{align*}
\Coder(S(M))&\nach\Der(R),\\
Q&\mapsto s^Q
\end{align*}
where $s^Q(g)=g\circ Q$,
is bijective and the restriction gives rise to an isomorphism
$$\Codiff(S(M))\nach\Diff(R).$$
\end{prop}
\begin{proof}
The injectivity is clear. Surjectivity:
A derivation $s$ of degree $j$ on $R$ induces a differential
(also denoted by $s$) on $\bR=S(M)^\ast$.
We have to find 
a coderivation $Q$ of degree $j$ on $S(M)$ such that,
for $u\in S(M)^\ast$, we have $s(u)=u\circ Q$.\\

For each $i\in I$, set $f_i:=s(x_i)$. Then,
$f_i$ is a product $((f_i)_n)_{n\geq 1}$ with 
$(f_i)_n\in\Hom^{-g(i)+1}(M^{\odot n},k)$.
We define the coderivation $Q$ by
$$Q_n(m_1,\ldots ,m_n):=\sum_{i\in I}(f_i)_n(m_1,\ldots ,m_n)\cdot e_i,$$
for homogeneous $m_1,\ldots ,m_n\in M$.
In fact, the non-vanishing terms in the sum satisfy the condition
$g(m_1)+\ldots +g(m_n)=g(i)$, hence the sum is finite.
To show that
for $u\in S(M)^\ast$, we have $s(u)=u\circ Q$, 
it is enough to show that for all $i\in I$,
$s(x_i)=x_i\circ Q$. But by definition, 
for $m_1,\ldots ,m_n\in M$,
we have
$$(x_i\circ Q)_n(m_1,\ldots ,m_n)=(f_i)_n(m_1,\ldots ,m_n)=(s(x_i))
(m_1,\ldots,m_n).$$
The second statement is a direct consequence of the first.
\end{proof}

As consequence, the differential $s$ on $R$ induces a codifferential
$Q^M$ on $S(M)$. We consider the pair $(M,Q^M)$ as formal
DG manifold in $\dgmanf^LG$. It has the following property:
The restriction of $Q^M$ to $M^0$ defines a formal
map $M^0\nach M^1$. Its zero locus is isomorphic to $X$.\\

Summarizing the above construction,
to each formal space germ $X$ with associated
formal analytic algebra $A$, we can construct a formal
DG manifold $(M,Q^M)$, containing $X$ as ``subspace''.
Of course, $(M,Q^M)$ depends on the choice of
the resolvent $(R,s)$. But we will show that 
$(M,Q^M)$ is well defined up to weak equivalence, i.e.
that the assignment $X\mapsto (M,Q^M)$ defines
a functor
$$F:\anf\nach\dgmanf^LG/\approx.$$

\index{$F$}

\begin{lemma}\label{mdgut1}
If $W=(W,d)$ is a DG $k$-vectorspace and if the dual complex\\
$\Hom(W,k)$ is acyclic, then $W$ is acyclic. Consequently,
if $f:V\nach W$ is a morphism of DG $k$-vectorspaces such that
the dual complex $f^\ast:W^\ast\nach V^\ast$ is a quasi-isomorphism,
then $f$ is a quasi-isomorphism.
\end{lemma}
\begin{proof}
Assume that $M$ is cyclic, i.e. there is an $n$
and an element $a\in M^n$ such that $d^n(a)=0$ and
$a\not\in\Im d^{n-1}$.
Let $B'$ be a base of $\im d^{n-1}$. We extend $B'\cup\{a\}$
to a base $B$ of $M^n$.
Let $p:M^n\nach k$ be the projection on the coordinate $a$ of $B$.
Then, $d^\ast(p)=p\circ d^{n-1}=0$ and $p(a)=1$,
hence $p\not\in\Im d^\ast$.
Contradiction !
\end{proof}

\begin{lemma}\label{mdgut2}
Let $f:M\nach M'$ be a morphism of formal DG-manifolds such
that the corresponding map $S(M)\nach S(M')$ is a quasi-isomorphism
of complexes. Then, $f$ is a weak equivalence.
\end{lemma}

\begin{proof}
By the Decomposition Theorem for $L_\infty$-algebras (see
Lemma 4.5 of \cite{Kont}), we may assume that $M$ is minimal and that
$f$ is strict. In this case, the homomorphism $f:S(M)\nach S(M')$
of DG coalgebras
is a direct sum of maps of complexes
$f_1:M\nach M'$ and 
$$\sum_{j\geq 2}f_1^{\odot j}:\coprod_{j\geq 2}M^{\odot j}\nach
\coprod_{j\geq 2}{M'}^{\odot j}.$$
Since the sum is a quasi-isomorphism, both factors
are quasi-isomorphisms.
\end{proof}

\begin{kor}
Let $F:(M,Q^M)\nach(M',Q^{M'})$ be a morphism of formal DG manifolds in
$\dgmanf^G$ and suppose that the dual map 
If $S(M')^\ast\nach S(M)^\ast$ is a quasi-isomorphism 
of free DG algebras,
then $F$ is a weak equivalence.
\end{kor}

\begin{proof}
This follows by Lemma~\ref{mdgut1}
and \ref{mdgut2}.
\end{proof}

Thus, we have proved the functoriality of $F$. Next, we 
define a functor\index{$V$}
$$V:\dgmanf^{GM}\nach\anf$$
as already mentioned above:
For a minimal DG manifold $(M,Q^M)$ in $\dgmanf^{MG}$, set $V(M,Q^M)$ to
be the zero locus of the formal map
$M^0\nach M^1$, induced by $Q^M$.
It can easily be seen that
the composition $V\circ F$ is the identity on $\anf$.
As a consequence, we get the following theorem:

\begin{satz}
The functor $F$ embeds
$\anf$ as full subcategory into $\dgmanf^{GM}$.
\end{satz}

\section{Deformations and embedded deformations}\label{embdef}

In this section we recall some classical results,
showing that each deformation of a singularity is 
equivalent to an embedded deformation.\\

A morphism $G:\sC\nach\sD$ of fibered gruppoids over
the category $\an$ of complex space germs is called
\textbf{smooth} if the following condition holds:
If $\beta:b\nach b'$ is a morphism in $\sD$ such that
$G(\beta):S\nach S'$ is a closed embedding, and if
$a$ is an object in $\sC$ such that $G(a)=b$,
then there is a morphism $\alpha:a\nach a'$ in $\sC$ such
that $G(\alpha)=\beta$.\\

Consider a complex space germ $X$ with
corresponding analytical algebra $\Oh_X$. Suppose that
$X$ is embedded in the smooth space germ $P$ with corresponding
analytic algebra $R^0$. Let $R=(R,s)$ be a g-finite, free
algebra resolution of $\Oh_X$ such that $R^0=\Oh_P$.\\

For any space germ $(S,\Oh_S)$, set $R_S:=R\hat{\ot}_{\CC}\Oh_S$
and
$$\sC(S):=\{\delta\in\Der^1(R_S,R_S)|\;\delta(0)=0
\text{ and }(s+\delta)^2=0\}$$
Furthermore, let $\sD(S)$ be the equivalence class of
deformations of $X$ with base $S$, i.e. the equivalence class
of all
flat morphisms $\sX\nach S$ such that there is a cartesian
diagram
\begin{equation}\label{deform}
\xymatrix{
X\ar[r]\ar[d] & \sX\ar[d]\\
\ast\ar[r] & S
}\end{equation}

Then, $\sC$ and $\sD$ are fibered gruppoids over $\an$ and
we define a morphism $G:\sC\nach\sD$ as follows:
For $\delta\in\sC(S)$, let $\sX$ be the space germ
with $\Oh_{\sX}=H^0(R_S,\delta+s)$ and $\sX\nach S$ the composition
of the closed embedding $\sX\nach S\times P$ and the canonical
projection
$S\times P\nach S$.
Obviously, there is a cartesian diagram (\ref{deform}).
I. e. $G(\delta):=\sX\nach S$ is a deformation of $X$.
We want to remind the proof of the well-known fact that $G$ is 
smooth.\\

Let $(A,\m)$ be a local analytic algebra, $B$ a graded, g-finite
free $A$-algebra and $C$ a flat DG-algebra over $A$. For $A$-modules
$M$, we set $M':=M\hat{\ot}_AA/\m$.
The following statement is a special case of
Proposition 8.20 in Chapter I of \cite{BiKo}: 
\begin{prop}\label{bick}
Let $v'\in\Der_{B'_0}^1(B',B')$ be a differential and
$\phi':B'\nach C'$ a surjective quasi-isomorphism of DG-algebras
over $A'$. Then, there is a differential $v\in\Der^1_{B_0}(B,B)$,
lifting $v'$ and a surjective quasi-isomorphism 
$\phi:B\nach C$ of DG-algebras, lifting
$\phi'$.
\end{prop}

\begin{kor}
For all $S$ in $\an$, $G(S):\sC(S)\nach\sD(S)$ is surjective.
\end{kor}

\begin{proof}
For $\sX\nach S$ in $\sD(S)$, we have to find
a $\Oh_S$-derivation $\delta:R_S\nach R_S$ of degree 
1 with $\delta(0)=0$ such that $\delta+s$ is a differential and
a surjective quasi-isomorphism
$(R_S,s+\delta)\nach \Oh_{\sX}$. Since $R_S\hat{\ot}_{\Oh_S}\CC=R$ and
$\Oh_{\sX}\hat{\ot}_{\Oh_S}\CC=\Oh_X$, the existence follows by
Proposition~\ref{bick}.
\end{proof}

\begin{kor}\label{smooth}
$G$ is smooth.
\end{kor}

\begin{proof}
We have to show that for each $\delta\in\sC(S)$ and
each morphism
$$\xymatrix{
\sX:=V(S\times P,\delta+s)\ar[r]\ar[d] & \sX'\ar[d]\\
S\ar[r] & S'
}$$
of deformations of $X$, there exist $\delta'\in\sC(S')$
such that $G(\delta')=\sX'$ and a cartesian
diagram
$$\xymatrix{
(R_{S'},\delta'+S)\ar[r] & (R_S,\delta+s)\\
\Oh_{S'}\ar[u]\ar[r] & \Oh_S\ar[u]
}$$
Setting $A:=\Oh_{S'}$, this follows by Proposition~\ref{bick}.
\end{proof}

In the literature (see \cite{BiKo}, for instance),
the deformation functor is defined such
that a space germ $S$ maps to the quotient of $\sC(S)$ by
the Lie group, associated to the Lie algebra
$\Der^0(R_S,R_S)$. In fact, $G$ factors through this
quotient and the first factor is even ``minimal
smooth''. For the construction here, we don't need to consider
this group action to get semi-universal deformations.
One can say that the group action is replaced by
the going - over to a minimal model.

\section{A formal semi-universal deformation}

In this section, we apply the new method for the construction
of a formal semi-universal deformations to 
isolated singularities $X$.
Let $(M,Q^M):=F(X)$ be the formal DG-manifold in $\dgmanf^MG$,
assigned to the space germ $X$. As in Section~\ref{functors}, denote
the resolvent of $A=\Oh_X$, having $S(M)^\ast$ as
completion, by $(R,s)$.\\

By Theorem 5.13 of \cite{Schuh1}, there is a semiuniversal deformation
$(V,Q^V,Q)$ of $(M,Q^M)$. Recall that as graded modules $V=H[1]$,
where $H$ denotes the cohomology of $\Coder(S(M),S(M))$, i.e. the
tangent cohomology of $X$. 
It is well-known that $H$ is g-finite.\\

We apply the functor $V$ to the
morphism $(V\times M,Q^V+Q^M+Q)\nach (V,Q^V)$
and get a morphism
$\sY\nach Y$ in $\anf$.

\begin{satz}
The morphism $\Y\nach Y$ is a formal semi-universal deformation
of the space germ $X$.
\end{satz}
\begin{proof}
Let
$$\xymatrix{
\sX\ar[d] & X\ar[l]\ar[d]\\
S & \ast\ar[l]
}$$ 
be any formal deformation of $X$. 
By Corollary~\ref{smooth}, there is a morphism of the deformation
$\sX\nach S$ to an embedded deformation $\tsX\nach S$,
where $\tsX$ is such that $\Oh_{\tsX}=H^0(R_S,s+\delta)$,
for a certain $\delta\in\sC(S)$ (see Section~\ref{embdef}).
I.e. there is a cartesian diagram
$$\xymatrix{
\tsX\ar[d] & \sX\ar[l]\ar[d]\\
S & S\ar[l]
}$$
Set $(B,Q^B):=F(S)$.
By Proposition~\ref{tan},
to $\delta$, there corresponds a coderivation
$Q_\delta$ in $\Coder^{+1}(S(B\times M),S(B\times M))$,
defining a deformation $(B,Q^B,Q_\delta)$ of $(M,Q^M)$.
Since $(V,Q^V)$ is semi-universal, there is a morphism
$$\xymatrix{
(B\times M,Q^B+Q^M+Q_\delta)\ar[r]\ar[d] & (V\times
M,Q^V+Q^M+Q)\ar[d]\\
(B,Q^B)\ar[r] & (V,Q^V)
}$$
of deformations. 
Application of the functor $V$ gives a cartesian diagram
$$\xymatrix{
\tsX\ar[r]\ar[d] & \sY\ar[d]\\
S\ar[r] & Y
}$$
which obviously respects the distinguished fiber $X\nach\ast$.
This shows that $\sY\nach Y$ is versal. Since
$Y$ is a formal analytic subgerm of $V^0=H^1$,
we have $\dim(TY)\leq \dim H^1$. Thus, necessarily $\sY\nach Y$
is semi-universal (see Chapter 2.6 of \cite{Pala}).
\end{proof}

\lz
\begin{center}
Institut Fourier\\
UMR 5582\\
BP 74\\
38402 Saint Martin d'H\`eres\\
France\\
\lz
frank.schuhmacher@ujf-grenoble.fr
\end{center}

\end{document}